\input amstex
\loadeufm

\documentstyle{amsppt}

\magnification=\magstep1

\baselineskip=20pt
\parskip=5.5pt
\hsize=6.5truein
\vsize=9truein
\NoBlackBoxes

\define\br{{\Bbb R}}

\define\e{{\varepsilon}}
\define\OO{{\Omega}}

\topmatter
\title
On Estimates of Biharmonic Functions on Lipschitz and Convex Domains
\endtitle

\author Zhongwei Shen
\endauthor

\leftheadtext{Zhongwei Shen}
\rightheadtext{Biharmonic Functions}
\address Department of Mathematics, University of Kentucky,
Lexington, KY 40506.
\endaddress

\email shenz\@ms.uky.edu
\endemail

\abstract
Using  Maz'ya  type integral identities with power weights, 
we obtain new boundary 
estimates for biharmonic functions
on Lipschitz and convex domains in $\br^n$.
For $n\ge 8$, combined with a result in \cite{S2}, these estimates
lead to the solvability of the $L^p$ Dirichlet problem
for the biharmonic equation on Lipschitz domains
for a new range of $p$. In the case of convex domains, the estimates allow
us to show that
the $L^p$ Dirichlet problem is uniquely solvable for
any $2-\e<p<\infty$ and $n\ge 4$.
\endabstract



\endtopmatter

\document

\centerline{\bf 1. Introduction}

Let $\OO$ be a bounded domain in $\br^n$ with Lipschitz boundary.
Let $N$ denote the outward unit normal to $\partial\OO$.
We consider the $L^p$ Dirichlet problem for the biharmonic
equation,
$$
\left\{
\aligned
&\Delta^2 u =0\ \ \ \ \text{ in } \OO,\\
&u=f \in L^p_1(\partial\OO),\ \frac{\partial u}{\partial N}=g
\in L^p(\partial\OO)\ \ \ \ \text{on } \partial\OO,\\
&(\nabla u)^* \in L^p(\partial\OO)
\endaligned
\right.
\tag 1.1
$$
where $L^p_1(\partial\OO)$ denotes the space of functions in
$L^p(\partial\OO)$ whose first order (tangential) derivatives are also in 
$L^p(\partial\OO)$. We point out that the boundary values in (1.1)
are taken in the sense of non-tangential convergence a.e.
with respect to the surface measure on $\partial\OO$. As such, one 
requires that the non-tangential maximal function $(\nabla u)^*$
is in $L^p(\partial\OO)$.

For $n\ge 2$, the Dirichlet problem (1.1) with $p=2$ was solved
by Dahlberg, Kenig and
Verchota \cite{DKV}, using bilinear estimates for harmonic
functions. The result was then extended
to the case $2-\e<p<2+\e$ by a real variable argument, where 
$\e>0$ depends on $n$ and $\OO$.
They also showed that the restriction $p>2-\e$ is necessary for
general Lipschitz domains.
In \cite{PV1,PV2}, Pipher and Verchota proved that if $n=3$ (or 2),
the $L^p$ Dirichlet problem (1.1) is uniquely solvable
for the sharp range $2-\e<p\le \infty$.
Moreover, they pointed out that (1.1) is {\it not} solvable in general
for $p>6$ if $n=4$, and for $p>4$ if $n\ge 5$.
Recently in \cite{S1,S2}, 
for $n\ge 4$ and $p$ in a certain range,
we established the
solvability of the $L^p$ Dirichlet problem
for higher order elliptic equations and systems,
using a new approach via $L^2$ estimates and weak
reverse H\"older inequalities.
In particular,
we were able to solve the $L^p$ Dirichlet problem (1.1) in the
following cases,
$$
\left\{
\alignedat3
 & 2-\e <p< 6+\e & \quad\quad\quad &\text{ for }\ \  n=4,\\
& 2-\e <p <4 +\e & \quad\quad\quad &\text{ for }\ \  n=5,6,7,\\
& 2-\e <p< 2+\frac{4}{n-3} +\e &\quad\quad\quad & \text{ for }
\ \  n\ge 8.
\endalignedat
\right.
\tag 1.2
$$
This gives the sharp ranges of $p$ for $4\le n\le 7$.
It should be pointed out that the sharp range $2-\e<p<4+\e$
for the case $n=6,7$ in (1.2) relies on a classical result
of Maz'ya \cite{M1,M2} on the boundary regularity of
biharmonic functions in arbitrary domains.
The approach we will use in this paper is inspired
by the work of Maz'ya \cite{M1,M2,M4}
(we shall come back to this point later).
We mention that if the domain $\OO$ is $C^1$, then (1.1) is 
uniquely solvable for all $n\ge 2$ and $1<p\le \infty$ \cite {CG,V1,PV2}.
For related work on the $L^p$ Dirichlet problem for
the polyharmonic equation
and general higher order equations and systems
on Lipschitz domains,
we refer the reader to \cite{V2,PV3,PV4,K,V3,S1,S2}.

The purpose of this paper is twofold.
First we study the case $n\ge 8$ for which the
question of the sharp ranges of $p$ remains open
for Lipschitz domains. Secondly we initiate the study of
the $L^p$ Dirichlet problem (1.1) on convex domains.
Note that any convex domain is Lipschitz, but may not be $C^1$.

Let $I(Q,r)=B(Q,r)\cap \partial\OO$ and $T(Q,r)=B(Q,r)\cap \OO$
where $Q\in \partial\OO$ and $r>0$. Our starting point is the following
theorem.

\proclaim{\bf Theorem 1.3}
Let $\OO$ be a bounded Lipschitz domain in $\br^n$, $n\ge 4$. Suppose that
there exist constants $C_0>0$, $R_0>0$ and $\lambda\in (0,n]$ such that
for any $0<r<R<R_0$ and $Q\in \partial\OO$,
$$
\int_{T(Q,r)} |\nabla v|^2\, dx
\le C_0\left(\frac{r}{R}\right)^\lambda
\int_{T(Q,R)}
|\nabla v|^2\, dx,
\tag 1.4
$$
whenever $v$ satisfies
$$
\left\{
\aligned
&\Delta^2 v=0\ \ \ \ \text{ in } \OO,\\
&v=\frac{\partial v}{\partial N} =0\ \ \ \text{ on } I(Q,R),\\
&(\nabla v)^* \in L^2(\partial\OO).
\endaligned
\right.
\tag 1.5
$$
Then the $L^p$ Dirichlet problem (1.1) is uniquely solvable for
$$
2<p<2+\frac{4}{n-\lambda}.
\tag 1.6
$$
 Moreover, the solution $u$ satisfies
$$
\| (\nabla u)^*\|_{L^p(\partial\OO)}
\le C\big\{
\|\nabla_t f\|_{L^p(\partial\OO)}
+\| g\|_{L^p(\partial\OO)}\big\},
\tag 1.7
$$
where $\nabla_t f$ denotes the tangential derivatives of $f$ on $\partial\OO$.
\endproclaim

Theorem 1.3 is a special case of Theorem 1.10 in \cite{S2}
for general higher order homogeneous elliptic equations and systems
with constant coefficients. It reduces the study of the $L^p$
Dirichlet problem to that of local $L^2$ estimates
near the boundary.
The main body of this paper will be devoted to such estimates.
In particular, we will prove that if $n\ge 8$, then estimate (1.4)
holds for some $\lambda>\lambda_n$, where
$$
\lambda_n=\frac{n+10 +2\sqrt{2(n^2-n+2)}}{7}.
\tag 1.8
$$
We will also show that
if $\OO$ is convex and $n\ge 4$, then (1.4) holds for any
$0<\lambda<n$.
Consequently, by Theorem 1.3, we obtain the following.

\proclaim{\bf Main Theorem}
Let $\OO$ be a bounded Lipschitz domain in $\br^n$.

\item{a)} If $n\ge 8$, the $L^p$ Dirichlet problem (1.1) is uniquely
solvable for 
$$
2-\e<p<2+\frac{4}{n-\lambda_n} +\e.
\tag 1.9
$$

\item{b)} If $n\ge 4$ and $\OO$ is convex, the $L^p$ Dirichlet problem (1.1)
is uniquely solvable for $2-\e<p<\infty$.
\endproclaim

We remark that in the case of Laplace's equation $\Delta u=0$, 
the Dirichlet problem
in $L^p$
is uniquely solvable on convex domains for
all $1<p\le \infty$.
This follows easily from the $L^\infty$ boundary estimates
on the first derivatives of the Green's functions.
Whether a similar result (the $L^\infty$ boundary estimate
on the second derivatives) holds for biharmonic functions remains open
for $n\ge 3$ (see \cite{KM} for the case $n=2$).
Note that part (b) of the Main Theorem as well as 
its proof gives the $C^\alpha$ boundary estimate
of $u$ for any $0<\alpha<1$.
This seems to be the first regularity
result for biharmonic functions
on general convex domains in $\br^n$, $n\ge 4$.

As we mentioned earlier, our approach to estimate (1.4) 
is motivated by the work of Maz'ya \cite{M1,M2,M4}.
It is based on certain
integral identities for
$$
\int_\OO \Delta^2 u\cdot u\, \frac{dx}{\rho^\alpha}
\ \ \ \ \text{ and }\ \ \ \ 
\int_{\OO} \Delta^2 u\cdot\frac{\partial u}{\partial \rho} \,
 \frac{dx}{\rho^{\alpha -1}},
\tag 1.10
$$
where $\rho=|x-Q|$ with $Q\in \partial\OO$ fixed.
See (2.13) and (3.1).
These identities with power weights allow us to control
the integrals
$$
\int_\OO |\nabla u|^2\, \frac{dx}{\rho^{\alpha +2}}
\ \ \ \ \text { and }\ \ \ \ 
\int_\OO |u|^2\, \frac{dx}{\rho^{\alpha +4}}
\tag 1.11
$$
for certain values of $\alpha$.
We point out that 
integral identity (2.13) with $\alpha=n-4$
appeared first in \cite{M1,M2},
where it was used to establish a Wiener's type
condition on the boundary continuity for the biharmonic equation
$\Delta^2 u=f$ on arbitrary domains in $\br^n$ for $n\le 7$.
Since the restriction $n\le 7$ in \cite{M1,M2}
is related to the positivity of 
a quadratic form (see (1.12) below), the idea to prove part
(a) of the Main Theorem is to use the identity
(2.13) for certain $\alpha<n-4$ in the case $n\ge 8$.
However, it should be pointed out that
the main novelty of this paper is 
the new identity (3.1), on which 
the proof of part (b) of the 
Main Theorem is based.
This identity allows us to estimate the integrals in (1.11) 
on convex domains for any
$\alpha<n-2$.
We remark that due to the lack of maximum
principles for higher order equations,
identitties such as (2.13) and (3.1)
are valuable tools in the study of boundary regularities
in nonsmooth domains.

Finally we mention that the results
 in \cite{M1,M2} were subsequently extended to the polyharmonic
equation \cite{MD,M3} and general higher order 
elliptic equations \cite{M4}.
Also, the related question of the
positivity of the quadratic form
$$
\int_{\br^n} (-\Delta)^\lambda u\cdot u\, \frac{dx}{|x|^{n-2\lambda}}
\ \ \text{ for all real
function }
u\in C_0^\infty(\br^n), 
\tag 1.12
$$
 has been studied systematically by Eilertsen \cite{E1,E2}
for all $\lambda\in (0,n/2)$.

\noindent{\bf Acknowledgment.}
The author is indebted to Jill Pipher for bring the paper \cite{M2}
to his attention, and for several helpful discussions.
The author also would like to thank Vladimir Maz'ya 
for pointing out the relevance of the papers
\cite{M1,E1,E2}.

\medskip

\centerline{\bf 2. Boundary Estimates on Lipschitz Domains}

The goal of this section is to prove part (a) of the Main Theorem.
We begin with a Cacciopoli's inequality. Recall that for
$Q\in \partial\OO$, $T(Q,R)
=B(Q,R)\cap \OO$ and $I(Q,R)=B(Q,R)\cap \partial\OO$.
We assume that $0<R<R_0$, where $R_0$ is a constant
depending on $\OO$ so that for any $Q\in \partial\OO$,
 $T(Q,4R_0)$ is given by the intersection
of $B(Q,4R_0)$ and the region above a Lipschitz graph,
after a possible rotation.

\proclaim{\bf Lemma 2.1}
Let $u\in W^{2,2}(T(Q,R))$ for some $Q\in \partial\OO$ and $0<R<R_0$.
Suppose that $\Delta^2 u=0$ in $T(Q,R)$ and 
$u=0$, $\nabla u=0$ on $I(Q,R)$. Then
$$
\frac{1}{r^2}\int_{T(Q,r)} |\nabla u|^2\, dx
+\int_{T(Q,r)} |\nabla^2 u|^2\, dx
\le \frac{C}{r^4}\int_{T(Q,2r)\setminus T(Q,r)}
|u|^2\, dx,
\tag 2.2
$$
where $0<r<R/4$.
\endproclaim

\demo{Proof}
Let $\eta$ be a smooth function on $\br^n$ such that
$\eta=1$ on $B(Q,r)$, supp$\, \eta\subset B(Q,2r)$ and
$|\nabla^k\eta|\le C/r^k$ for $0\le k\le 4$.
Since $u\in W^{2,2}(T(Q,R))$ and $u=0$, $\nabla u=0$
on $I(Q,R)$, we have $u\eta^2\in W^{2,2}_0(\OO)$.
We will show that for any $\e>0$,
$$
\aligned
\int_\OO |\nabla^2(u\eta^2)|^2\, dx
\le & \e\int_\OO |\nabla^2(u\eta^2 )|^2\, dx
+\frac{\e}{r^2}
\int_\OO |\nabla (u\eta^2)|^2 \, dx\\
&\ \ \ \ \ \ \ \ \ 
+\frac{C_\e}{r^4} \int_{T(Q,2r)\setminus T(Q,r)} |u|^2\, dx.
\endaligned
\tag 2.3
$$
This, together with the Poincar\'e inequality
$$
\int_{T(Q,2r)} |\nabla(u\eta^2)|^2\, dx
\le C\, r^2\int_{T(Q,2r)} |\nabla^2(u\eta^2)|^2\, dx,
\tag 2.4
$$
yields the estimate (2.2). 

To prove (2.3), we use integration by parts and $\Delta^2 u=0$
in $T(Q,2r)$ to obtain
$$
\int_\OO|\nabla^2(u\eta^2)|^2\, dx
=\int_\OO |\Delta(u\eta^2)|^2\, dx
=\int_\OO \big\{ |\Delta(u\eta^2)|^2
-\Delta u\cdot \Delta (u\eta^4)\big\}\, dx.
\tag 2.5
$$
A direct computation shows that
$$
\aligned
&\Delta(u\eta^2)\cdot \Delta(u\eta^2)-\Delta u\cdot \Delta(u\eta^4)\\
&=u\Delta(u\eta^2)\Delta\eta^2 +4|\nabla u\cdot \nabla \eta^2|^2
+2u(\nabla u\cdot \nabla\eta^2)\Delta \eta^2
-u\Delta u(2|\nabla\eta^2|^2 +\eta^2\Delta\eta^2),
\endaligned
\tag 2.6
$$
In view of (2.3), the integral of the first term in the right side of (2.6)
can be handled easily by H\"older's inequality with an $\e$.
The remaining terms may be handled by using
 integration by parts, together with the following
observation.
For terms with 
$u\frac{\partial u}{\partial x_i} $, like the third term
in the right side of (2.6),  we may write 
$$
u\frac{\partial u}{\partial x_i}\psi
=\frac{1}{2} \frac{\partial }{\partial x_i}
\left(|u|^2\psi\right) -\frac{1}{2} |u|^2 \frac{\partial \psi}
{\partial x_i}.
\tag 2.7
$$
For terms with $\eta^2\frac{\partial u}{\partial x_i}\frac{\partial u}{
\partial x_j}$, like the second term, we use
$$
\aligned
\frac{\partial u}{\partial x_i}\frac{\partial u}{
\partial x_j}\eta^2\psi
& =\frac{\partial}{\partial x_i}
\left(\frac{\partial (u\eta^2)}{\partial x_j}\cdot u\psi\right)
-\frac{\partial^2(u\eta^2)}{\partial x_i\partial x_j}
\cdot u\psi
-u\frac{\partial u}{\partial x_j} \frac{\partial \psi}{\partial x_i}
\eta^2\\
&\ \ \ \ \ \ 
-u\frac{\partial u}{\partial x_i} \frac{\partial \eta^2}{\partial x_j}
\psi
-u^2\frac{\partial\eta^2}{\partial x_j}\frac{\partial \psi}{
\partial x_i}.
\endaligned
\tag 2.8
$$
Finally, for the last term which contains $\eta^2 u\Delta u$, we note that
$$
\eta^2 u\Delta u \cdot \psi
=\frac{\partial }{\partial x_i}
\left(\eta^2 u\frac{\partial u}{\partial x_i}\psi\right)
-u\frac{\partial u}{\partial x_i}\frac{\partial (\eta^2\psi)}
{\partial x_i}
-\eta^2 |\nabla u|^2 \psi.
\tag 2.9
$$
The rest of the proof, which we omit, is
fairly straightforward.
\enddemo

\remark{\bf Remark 2.10}
It follows from Lemma 2.1 that for any $0<r<R/2$ and $\alpha\in\br$,
$$
\int_{T(Q,r)} \frac{|\nabla u(x)|^2}{|x-Q|^{\alpha +2}}\, dx
+\int_{T(Q,r)}\frac{|\nabla^2 u(x)|^2}{|x-Q|^\alpha}
\, dx
\le C\int_{T(Q,2r)}
\frac{|u(x)|^2\, dx}{|x-Q|^{\alpha +4}}.
\tag 2.11
$$
This may be seen by writing $T(Q,r)$ as $\cup_{j=0}^\infty
T(Q, 2^{-j}r)\setminus T(Q,2^{-j-1} r)$.
\endremark

The key step to establish estimate (1.4) relies on
the following extension of an integral identity due to Maz'ya \cite{M1,M2}.
 
\proclaim{\bf Lemma 2.12}
Suppose that $u\in C^2(\overline{\OO})$ and 
$u=0$, $\nabla u=0$ on $\partial\OO$.
Then for any $\alpha\in \br$,
$$
\aligned
\int_\OO \Delta u\cdot \Delta \big(\frac{u}{\rho^\alpha}\big)\, dx
&=\int_\OO |\Delta u|^2\, \frac{dx}{\rho^\alpha}
+2\alpha \int_\OO |\nabla u|^2 \, \frac{dx}{\rho^{\alpha +2}}
-2\alpha(\alpha +2)\int_\OO \big|\frac{\partial u}{\partial
\rho}\big|^2\, \frac{dx}{\rho^{\alpha +2}}\\
&\ \ \ \ \ +
\frac{1}{2}\alpha (\alpha +2)(n-2-\alpha)
(n-4-\alpha)
\int_\OO |u|^2\, \frac{dx}{\rho^{\alpha +4}},
\endaligned
\tag 2.13
$$
where $\rho=|x-y|$ and $\frac{\partial u}{\partial\rho}=<\nabla u(x),
(x-y)/\rho>$ with $y\in \overline{\OO}^c$ fixed.
\endproclaim

\demo{Proof} 
We will use the summation convention that the repeated indices
are summed from $1$ to $n$.
First, note that
$$
\aligned
\int_\OO \Delta u\cdot \Delta\big(\frac{u}{\rho^\alpha}\big)
dx
&=\int_\OO |\Delta u|^2\, \frac{dx}{\rho^\alpha}
+2\int_\OO \Delta u \cdot \frac{\partial u}{\partial x_j}
\cdot \frac{\partial }{\partial x_j} \big(\frac{1}{\rho^\alpha}\big)\, dx\\
&\ \ \ \ \ \
+\int_\OO u\Delta u \cdot \Delta \big(\frac{1}{\rho^\alpha}\big)\, dx.
\endaligned
\tag 2.14
$$
Next it follows from integration by parts that
$$
2\int_\OO \Delta u\cdot\frac{\partial u}{\partial x_j}
\frac{\partial }{\partial x_j}\big(\frac{1}{\rho^\alpha}\big)\, dx
=\int_\OO |\nabla u|^2\Delta\big(\frac{1}{\rho^\alpha}\big)\, dx
-2\int_\OO \frac{\partial u}{\partial x_i}
\frac{\partial u}{\partial x_j}
\frac{\partial^2}{\partial x_i\partial x_j}
\big(\frac{1}{\rho^\alpha}\big)\, dx.
\tag 2.15
$$
Similarly, we have
$$
\int_\OO u\Delta u \cdot\Delta\big(\frac{1}{\rho^\alpha}\big)\, dx
=-\int_\OO |\nabla u|^2\, \Delta\big(\frac{1}{\rho^\alpha}\big)\, dx
+\frac12 \int_\OO
|u|^2 \, \Delta^2\big(\frac{1}{\rho^\alpha}\big)\, dx.
\tag 2.16
$$
Substituting (2.15) and (2.16) into (2.14), we obtain
$$
\aligned
\int_\OO \Delta u\cdot \Delta \big(\frac{u}{\rho^\alpha}\big)\, dx
&=\int_\OO |\Delta u|^2\, \frac{dx}{\rho^\alpha}
-2\int_\OO \frac{\partial u}{\partial x_i}
\cdot\frac{\partial u}{\partial x_j}\cdot\frac{\partial^2}{\partial x_i
\partial x_j}\big(\frac{1}{\rho^\alpha}\big)\, dx\\
&\ \ \ \ \
+\frac12 \int_\OO |u|^2\, \Delta^2 \big(\frac{1}{\rho^\alpha}\big)\, dx.
\endaligned
$$
The desired formula (2.13) now follows from the fact that
$$
\aligned
&\frac{\partial^2}{\partial x_i\partial x_j}
\big(\frac{1}{\rho^\alpha}\big)
=-\alpha \rho^{-\alpha -2} \delta_{ij}
+\alpha (\alpha+2) (x_i-y_i)(x_j-y_j)\rho^{-\alpha -4},\\
& \Delta^2 \big(\frac{1}{\rho^\alpha}\big)
=\alpha (\alpha +2)(n-2-\alpha)
(n-4-\alpha )\rho^{-\alpha -4},
\endaligned
\tag 2.17
$$
for any $\rho=|x-y|\neq 0$. The proof is complete.
\enddemo

\proclaim{\bf Lemma 2.18} Under the same assumption as in Lemma 2.12,
we have
$$
\int_\OO \Delta u\cdot\frac{\partial u}{\partial\rho}
\frac{dx}{\rho^{\alpha +1}}
=\frac{1}{2}(n-4-\alpha)
\int_\OO |\nabla u|^2\,\frac{dx}{\rho^{\alpha +2}}
+(\alpha +2)\int_\OO \big|\frac{\partial u}{\partial\rho}\big|^2\,
\frac{dx}{\rho^{\alpha +2}}.
\tag 2.19
 $$
\endproclaim

\demo{Proof}
It follows from integration by parts that
$$
\aligned
\int_\OO \Delta u\cdot\frac{\partial u}{\partial\rho}
\frac{dx}{\rho^{\alpha +1}}
&=\int_\OO \Delta u\cdot \frac{\partial u}{\partial x_i}\,
(x_i-y_i)\frac{dx}{\rho^{\alpha +2}}\\
&=\frac12 \int_\OO
|\nabla u|^2\, \frac{\partial }{\partial x_i}
\big(\frac{x_i-y_i}{\rho^{\alpha +2}}\big)\, dx
-\int_\OO
\frac{\partial u}{\partial x_i}\cdot \frac{\partial u}{\partial x_j}
\cdot \frac{\partial }{\partial x_j}
\big(\frac{x_i-y_i}{\rho^{\alpha +2}}\big)\, dx\\
&=\frac{1}{2}(n-4-\alpha)
\int_\OO |\nabla u|^2\,\frac{dx}{\rho^{\alpha +2}}
+(\alpha +2)\int_\OO \big|\frac{\partial u}{\partial\rho}\big|^2\,
\frac{dx}{\rho^{\alpha +2}}.
\endaligned
$$
\enddemo

Lemma 2.12, together with Lemma 2.18,  allows us to estimate
$$
\int_\OO |u|^2\, \frac{dx}{\rho^{\alpha+4}}
\ \ \text{ and }\ \ 
\int_{\OO} |\nabla u|^2\, \frac{dx}{\rho^{\alpha +2}}
\ \ \ \ \text{ by }\ \ \
\int_\OO \Delta u\cdot \Delta \big(\frac{u}{\rho^\alpha}
\big)\, dx
$$
for certain values of $\alpha$. 

\proclaim{\bf Lemma 2.20} Let $\OO$ be a bounded Lipschitz domain
in $\br^n$, $n\ge 5$. Suppose that $u\in C^2(\overline{\OO})$ and
$u=0$, $\nabla u=0$ on $\partial\OO$. 
Then, if $0<\alpha\le n-4$ and $n^2+2n\alpha -7\alpha^2 -8\alpha>0$,
we
have
$$
\int_\OO |\nabla u|^2 \, \frac{dx}{\rho^{\alpha +2}}
\le C_{n,\alpha}
\int_\OO \Delta u \cdot \Delta\big(\frac{u}{\rho^\alpha}\big)\, dx,
\tag 2.21
$$
where $C_{n,\alpha}>0$ depends only on $n$ and $\alpha$.
\endproclaim

\demo{Proof} We first use (2.19) for $0<\alpha\le n-4$ to obtain
$$
\aligned
\frac{n+\alpha}{2}\int_\OO \big|\frac{\partial u}{\partial \rho}\big|^2\,
\frac{dx}{\rho^{\alpha +2}}
&\le \int_{\OO} \Delta u\cdot \frac{\partial u}{\partial\rho}\,
\frac{dx}{\rho^{\alpha +1}}\\
&\le   
\left\{\int_\OO |\Delta u|^2\, \frac{dx}{\rho^\alpha}\right\}^{1/2}
\left\{ \int_\OO |\frac{\partial u}{\partial\rho}|^2
\,\frac{dx}{\rho^{\alpha +2}}\right\}^{1/2},
\endaligned
$$
where the Cauchy inequality is also used. It follows that
$$
\frac14 (n+\alpha )^2\int_\OO \big|\frac{\partial u}{\partial
\rho}\big|^2\, \frac{dx}{\rho^{\alpha +2}}
\le \int_\OO |\Delta u|^2 \, \frac{dx}{\rho^\alpha}.
\tag 2.22
$$
Since $0<\alpha\le n-4$, in view of (2.13) and (2.22), we have
$$
\aligned
\int_\OO \Delta u\cdot \Delta \big(\frac{u}{\rho^\alpha}\big)\,
dx
&\ge \left\{ \frac14 (n+\alpha)^2 +2\alpha -2\alpha (\alpha +2)\right\}
\int_\OO \big|\frac{\partial u}{\partial \rho}\big|^2\, 
\frac{dx}{\rho^{\alpha +2}}\\
&
=\frac{1}{4}
(n^2 +2n\alpha -7\alpha^2 -8\alpha)
\int_\OO \big|\frac{\partial u}{\partial \rho} \big|^2\, 
\frac{dx}{\rho^{\alpha +2}}.
\endaligned
 $$
Thus, if $n^2 +2n \alpha -7\alpha^2 -8\alpha>0$, by (2.13) again,
$$
\aligned
2\alpha\int_\OO |\nabla u|^2\, \frac{dx}{\rho^{\alpha+2}}
&\le \int_\OO \Delta u\cdot \Delta \big(\frac{u}{\rho^\alpha}
\big)\, dx
+2\alpha (\alpha +2)\int_\OO
\big|\frac{\partial u}{\partial \rho}\big|^2\, \frac{dx}{\rho^{\alpha +2}}\\
&\le C\, \int_\OO \Delta u\cdot \Delta \big(\frac{u}{\rho^\alpha}
\big)\, dx.
\endaligned
$$
The proof is finished.

\remark{\bf Remark 2.23}
Let $\alpha=n-4$. Then
$n^2+2n\alpha-7\alpha^2-8\alpha= 4(-n^2+10n -20)>0$
for $n=5,6,7$. It follows that (2.21) holds for 
$\alpha=n-4$ in the case $n=5$, $6$ or $7$.
This was the result obtained by Maz'ya in \cite{M1,M2}. If $n\ge 8$, then (2.21)
holds for $0<\alpha<\alpha_n<n-4$, where
$$
\alpha_n=\frac17 (n-4+2\sqrt{2(n^2-n+2)})
\tag 2.24
$$
is the positive root of $n^2+2n\alpha -7\alpha^2 -8\alpha=0$.
\endremark

\remark{\bf Remark 2.25}
If $n\ge 8$ and $\alpha=\alpha_n$ given by (2.24), we observe that
the first three terms on the right side of (2.13) is
nonnegative, by an inspection of th proof of Lemma 2.20.
It follows that
$$
\int_\OO |u|^2\, \frac{dx}{\rho^{\alpha +4}}
\le C_n\int_\OO \Delta u\cdot \Delta \big(\frac{u}{
\rho^\alpha}\big)\, dx.
\tag 2.26
$$
Since $C_0^\infty(\OO)$ is dense in $W^{2,2}_0(\OO)$,
inequality (2.26) holds for any $u\in W^{2,2}_0(\OO)$.
\endremark

We are now in a position to give the proof of part (a) of the Main Theorem.

\proclaim{\bf Theorem 2.27}
Let $\OO$ be a bounded Lipschitz domain in $\br^n$, $n\ge 8$. Then the
$L^p$ Dirichlet problem (1.1) is uniquely solvable for
$2-\e<p<2+\frac{4}{n-\lambda_n} +\e$,
where $\lambda_n=\alpha_n+2$ is given in (1.8).
\endproclaim

\demo{Proof}
By Theorem 1.3, we only need to show that estimate (1.4) holds for some
$\lambda>\lambda_n=\alpha_n +2$. To this end, we fix $Q\in\partial\OO$
and $0<R<R_0$, where $R_0$ is a constant depending on $\OO$.
Let $v$ be a function on $\OO$ satisfying (1.5).
Let $\eta$ be a smooth function on $\br^n$ such that
$\eta=1$ on $B(Q,r)$, supp$\, \eta\subset B(Q,2r)$ and
$|\nabla^k\eta|\le C/r^k$ for $0\le k\le 4$ where
$0<r<R/4$. Since $v=\frac{\partial v}{\partial N}=0$
on $I(Q,R)$ and $(\nabla v)^*\in L^2(\partial\OO)$, by the
regularity estimate $\|(\nabla\nabla v)^*\|_2
\le C\, \|\nabla_t\nabla v\|_2$ established in \cite{V2}, we know
$v\eta\in W^{2,2}_0(\OO)$. Thus we may apply estimate
(2.26) to $u=v\eta$ with $\alpha=\alpha_n$ and
$\rho=|x-y|$, where $y\in \overline{\OO}^c$. We obtain
$$
\int_\OO |v\eta|^2\, \frac{dx}{\rho^{\alpha +4}}
\le C\, \int_\OO 
\Delta (v\eta)\cdot \Delta \big(v\eta\rho^{-\alpha}\big)\, dx.
\tag 2.28
$$
Using an identity similar to (2.6), 
$$
\aligned
&\Delta(v\eta)\cdot \Delta(v\eta\rho^{-\alpha})
-\Delta v\cdot \Delta (v\eta^2 \rho^{-\alpha})\\
&=
v\Delta (v\eta \rho^{-\alpha})\Delta\eta
+2(\nabla v\cdot\nabla \eta)
\Delta(v\eta \rho^{-\alpha})
-2\big(\nabla(v\eta\rho^{-\alpha})\cdot\nabla \eta\big)
\Delta v
-v\eta\rho^{-\alpha} \Delta v\cdot \Delta \eta,
\endaligned
$$
and $\Delta^2 v=0$ in $\OO$, we get
$$
\aligned
\int_\OO |v\eta|^2\, \frac{dx}{\rho^{\alpha +4}}
\le C\, \int_\OO
&\bigg\{
v\Delta (v\eta \rho^{-\alpha})\Delta\eta
+2(\nabla v\cdot\nabla \eta)
\Delta(v\eta \rho^{-\alpha})\\
&\ \ \ \ -2\big(\nabla(v\eta\rho^{-\alpha})\cdot\nabla \eta\big)
\Delta v
-v\eta\rho^{-\alpha} \Delta v\cdot \Delta \eta\bigg\}\,dx.
\endaligned
\tag 2.29
$$
Note that $\rho^{-\alpha}$ and its derivatives are uniformly bounded for
$y\in B(Q,r/2)\setminus\overline{\OO}$ and
$x\in \text{supp}\,(|\nabla\eta|)\subset \{ x\in \br^n: r\le 
|x-Q|\le 2r\}$. It follows by a simple limiting argument that
(2.29) holds for $\rho=|x-Q|$. This gives
$$
\aligned
\int_{T(Q,r)} \frac{|v(x)|^2\,
dx}{|x-Q|^{\alpha +4}}
&\le \frac{C}{r^{\alpha +4}}
\int_{T(Q,2r)}
\left\{ |v|^2 +r^2 |\nabla v|^2
+r^4|\nabla^2 v|^2\right\}\, dx\\
&\le \frac{C}{r^{\alpha +4}}
\int_{T(Q,4r)\setminus T(Q,2r)}
|v|^2\, dx\\
&\le C_1\int_{T(Q,4r)\setminus T(Q,r)} \frac{|v(x)|^2\, dx}{
|x-Q|^{\alpha +4}},
\endaligned
\tag 2.30
$$
where the second inequality follows from 
Cacciopoli's inequality (2.2).
By ``filling'' the hole in (2.30), we obtain
$$
\int_{T(Q,r)} \frac{|v(x)|^2\, dx}
{|x-Q|^{\alpha +4}}
\le \frac{C_1}{C_1+1}
\int_{T(Q,4r)} \frac{|v(x)|^2\, dx}
{|x-Q|^{\alpha +4}}.
$$
This implies that there exists $\delta>0$ such that
$$
\aligned
\int_{T(Q,r)} \frac{|v(x)|^2\, dx}
{|x-Q|^{\alpha +4}}
&\le C\left(\frac{r}{R}\right)^\delta
\int_{T(Q,R/4)} \frac{|v(x)|^2\, dx}
{|x-Q|^{\alpha +4}}\\
&\le C\left(\frac{r}{R}\right)^\delta
\cdot\frac{1}{R^{\alpha +4}}
\int_{T(Q,R)} |v(x)|^2\, dx,     
\endaligned
$$
for any $0<r<R/4$, where the second inequality follows from (2.30).
Consequently,
$$
\int_{T(Q,r)} |v(x)|^2\,dx
\le C\left(\frac{r}{R}\right)^{\alpha +4 +\delta}
\int_{T(Q,R)} |v(x)|^2\, dx.
$$
This, together with
 Cacciopoli's inequality and Poincar\'e inequality, gives
$$
\aligned
\int_{T(Q,r/2)} |\nabla v|^2\, dx
&\le \frac{C}{r^2}\int_{T(Q,r)} |v|^2\, dx
\le \frac{C}{r^2} \left(\frac{r}{R}\right)^{\alpha +4+\delta}
\int_{T(Q,R)} |v|^2\, dx\\
&\le C\left(\frac{r}{R}\right)^{\alpha +2+\delta}
\int_{T(Q,R)} |\nabla v|^2\, dx.
\endaligned
$$
Thus we have established estimate (1.4) for $\lambda=\alpha_n+2+\delta
=\lambda_n +\delta$.
The proof is finished.
\enddemo  

\medskip

\centerline{\bf 3. Boundary Estimates on Convex Domains}

In this section we give the proof of part (b) of the Main Theorem.
By Theorem 1.3, it suffices to show that estimate (1.4)
holds for any $\lambda<n$. To do this,
the crucial step is to establish the following new integral identity,
$$
\aligned
&(\alpha +4-n)\int_\OO \Delta u\cdot\Delta\bigg(\frac{u}{\rho^\alpha}
\bigg)\, dx
-2\int_\OO \Delta^2 u\cdot \frac{\partial u}{\partial \rho}\,
\frac{dx}{\rho^{\alpha-1}}\\
&=\int_{\partial\OO} |\nabla^2 u|^2 <\frac{x-y}{\rho}, N>\,
\frac{d\sigma}{\rho^{\alpha-1}}
+4\alpha \int_\OO \big|\frac{\partial}{\partial\rho}
\left(\rho^{\frac{n-\alpha-2}{2}}\nabla u\right)\big|^2\, 
\frac{dx}{\rho^{n-2}}\\
&\ \ \ \ \ \ \ \ \ \ +2\alpha (\alpha +2) (n-\alpha-2)
\int_\OO \big|\frac{\partial}{\partial\rho}
\left(\rho^{\frac{n-\alpha-4}{2}} u\right)\big|^2\, \frac{dx}{\rho^{n-2}},
\endaligned
\tag 3.1
$$
where $u\in C^4(\overline{\OO})$ and
$u=0$, $\nabla u=0$ on $\partial\OO$.
Recall that $N$ denotes the outward unit normal to $\partial\OO$.
Also in (3.1), as before, $\rho=\rho(x)=
|x-y|$, $\frac{\partial u}{\partial\rho}
=<\nabla u,(x-y)/\rho>$ with $y\in \overline{\OO}^c$ fixed.
By a limiting argument, it is not hard to see that if $\alpha<n$,
(3.1) holds also for $y\in \partial\OO$. We will use (3.1)
with $\alpha =n-2$ for convex domain 
$\OO$. The key observation is that if $\OO$ is convex,
the boundary integral in (3.1) is nonnegative.
This is because $<P-Q,N(P)>\ge 0$ for any $P,Q\in\partial\OO$.
 
The proof of (3.1), which involves the repeated use of integration
by parts, will be given through a series of lemmas.

\proclaim{\bf Lemma 3.2}
Suppose $u\in C^2(\overline{\OO})$ and $u=0$, $\nabla u=0$ on $\partial\OO$.
Then, for any $\alpha\in \br$,
$$
\aligned
&\int_\OO \frac{\partial^2 u}{\partial x_i\partial x_j}
\cdot\frac{\partial^2}{\partial x_i\partial x_j}
\left(\frac{u}{\rho^\alpha}\right)\, dx
=\int_\OO |\nabla^2 u|^2\, \frac{dx}{\rho^\alpha}
+\alpha (n-\alpha-1)
\int_\OO |\nabla u|^2\, \frac{dx}{\rho^{\alpha +2}}\\
& \ \ \ 
-\alpha(\alpha +2)\int_\OO \big|\frac{\partial u}{\partial \rho}\big|^2
\, \frac{dx}{\rho^{\alpha +2}}
+\frac12 \alpha (\alpha +2)(n-\alpha -2)
(n-\alpha -4)\int_\OO |u|^2\, \frac{dx}{\rho^{\alpha +4}},
\endaligned
\tag 3.3
$$
where the repeated indices are summed from $1$ to $n$.
\endproclaim

\demo{Proof}
First we note that
$$
\aligned
& \int_\OO \frac{\partial^2 u}{\partial x_i\partial x_j}
\cdot\frac{\partial^2}{\partial x_i\partial x_j}
\left(\frac{u}{\rho^\alpha}\right)\, dx\\
&=\int_\OO \frac{\partial^2 u}{\partial x_i\partial x_j}
\cdot \bigg\{
\frac{\partial^2 u}{\partial x_i\partial x_j}\, \frac{1}{\rho^\alpha}
+2\,\frac{\partial u}{\partial x_i}\cdot
\frac{\partial\rho^{-\alpha}}{\partial x_j} 
+u\, \frac{\partial^2\rho^{-\alpha}
}{\partial x_i\partial x_j}\bigg\}\, dx.
\endaligned
\tag 3.4
$$
Next it follows from integration by parts  and 
$u=0$, $\nabla u=0$ on $\partial\OO$
that
$$
2\int_\OO \frac{\partial^2 u}{\partial x_i\partial x_j}
\cdot \frac{\partial u}{\partial x_i}
\cdot \frac{\partial \rho^{-\alpha}}{\partial x_j}\, dx
=-\int_\OO |\nabla u|^2 \, \Delta(\rho^{-\alpha})\, dx.
\tag 3.5
$$
and
$$
\aligned
&\int_\OO \frac{\partial^2 u}{\partial x_i\partial x_j}
\cdot u\cdot \frac{\partial^2\rho^{-\alpha}}
{\partial x_i\partial x_j}\, dx\\
&
=-\int_\OO \frac{\partial u}{\partial x_i}\cdot \frac{\partial u}
{\partial x_j} \cdot \frac{\partial^2 \rho^{-\alpha}}
{\partial x_i\partial x_j}\, dx
+\frac12 \int_\OO |u|^2\, \Delta^2 (\rho^{-\alpha})\, dx.
\endaligned
\tag 3.6
$$
Substituting (3.5) and (3.6) into (3.4), we obtain
$$
\aligned
&\int_\OO \frac{\partial^2 u}{\partial x_i\partial x_j}
\cdot\frac{\partial^2}{\partial x_i\partial x_j}
\left(\frac{u}{\rho^\alpha}\right)\, dx
=\int_\OO |\nabla^2 u|^2\, \frac{dx}{\rho^\alpha}
-\int_\OO |\nabla u|^2\, \Delta(\rho^{-\alpha})\, dx\\
&\ \ \ \ \ \ \ \ \ \ \ \ \ \ \ \ \ \
-\int_\OO \frac{\partial u}{\partial x_i}\,
\frac{\partial u}{\partial x_j}
\cdot\frac{\partial^2 \rho^{-\alpha}}
{\partial x_i\partial x_j}\, dx
+\frac12
\int_\OO |u|^2\, \Delta^2 (\rho^{-\alpha})\, dx.
\endaligned
$$
The desired formula now follows from this and (2.17).
\enddemo

\proclaim{\bf Lemma 3.7}
Suppose $u\in C^4(\overline{\OO})$ and $u=0$, $\nabla u=0$ on $
\partial\OO$. Then, for any $\alpha\in \br$,
$$
\aligned
&\int_\OO \Delta^2 u\cdot \frac{\partial u}{\partial \rho}
\, \frac{dx}{\rho^{\alpha-1}}
=-\frac12 \int_{\partial\OO}
|\nabla^2 u|^2 <\frac{x-y}{\rho}, N>\, \frac{d\sigma}{\rho^{\alpha-1}}\\
&\ \ \ \ \ +\frac12 (\alpha +4-n)\int_\OO |\nabla^2 u|^2\, \frac{dx}
{\rho^\alpha}
-2\alpha \int_\OO \big|\frac{\partial }{\partial\rho}
\nabla u\big|^2\, \frac{dx}{\rho^\alpha}\\
&\ \ \ \ \ \ +\frac12 \alpha (n-\alpha)\int_\OO |\nabla u|^2\, \frac{dx}{
\rho^{\alpha +2}}
-\frac12 \alpha (\alpha +2) (n-\alpha) 
\int_\OO \big|\frac{\partial u}{\partial \rho}\big|^2\,
\frac{dx}{\rho^{\alpha +2}},
\endaligned
\tag 3.8
$$
where $\rho=|x-y|$ with $y\in \overline{\OO}^c$ fixed
\endproclaim

\demo{Proof} By translation we may assume that $y=0$. 
Using integration by parts, we obtain
$$
\aligned
&\int_\OO \Delta^2 u \cdot \frac{\partial u}{\partial\rho}
\, \frac{dx}{\rho^{\alpha-1}}
=\int_\OO \frac{\partial^4 u}{\partial x_i^2
\partial x_j^2}\cdot\frac{\partial u}{\partial x_k}
\cdot\frac{x_k}{\rho^\alpha}\, dx\\
&
=-\int_\OO \frac{\partial^3 u}{\partial x_i\partial x_j^2}
\cdot \frac{\partial^2 u}{\partial x_i \partial x_k}\cdot
\frac{x_k}{\rho^\alpha}\, dx
-\int_\OO
\frac{\partial^3}{\partial x_i\partial x_j^2}\cdot
\frac{\partial u}{\partial x_k}\cdot
\frac{\partial }{\partial x_i}
\left(\frac{x_k}{\rho^\alpha}\right)\, dx.
\endaligned
\tag 3.9
$$
For the first term on the right side of (3.9),
again from integration by parts, 
we have
$$
\aligned
&-\int_\OO \frac{\partial^3 u}{\partial x_i\partial x_j^2}
\cdot \frac{\partial^2u}{\partial x_i\partial x_k}\cdot \frac
{x_k}{\rho^\alpha}\, dx
=-\frac12 \int_{\partial \OO}
|\nabla^2 u|^2\, <\frac{x}{\rho}, N>\, \frac{d\sigma}{\rho^{\alpha-1}}\\
&\ \ \ \ \ \ \ \ \ \
-\frac12 \int_\OO |\nabla^2 u|^2\, \frac{\partial }{\partial x_k}
\left(\frac{x_k}{\rho^\alpha}\right)\, dx
+
\int_\OO \frac{\partial^2 u}{\partial x_i\partial x_j}\cdot
\frac{\partial^2 u}{\partial x_i\partial x_k}\cdot
\frac{\partial }{\partial x_j}\left(\frac{x_k}{\rho^\alpha}
\right)\, dx,
\endaligned
\tag 3.10
$$
where we also used the observation that $\nabla u=0$ on $\partial\OO$
implies
$$
\frac{\partial }{\partial N}\left(\frac{\partial u}
{\partial x_i}\right)\cdot \frac{\partial }{\partial \rho}
\left(\frac{\partial u}{\partial x_i}\right)
= |\nabla^2 u|^2\, <\frac{x}{\rho},N>.
\tag 3.11
$$
For the second term on the right side of (3.9), we have
$$
\aligned
&-\int_\OO \frac{\partial^3 u}{\partial x_i\partial x_j^2}
\cdot\frac{\partial u}{\partial x_k} \cdot
\frac{\partial }{\partial x_i}
\left(\frac{x_k}{\rho^\alpha}\right)\, dx\\
&=\int_\OO \frac{\partial^2 u}{\partial x_i\partial x_j}
\cdot \frac{\partial^2 u}{\partial x_k\partial x_j}
\cdot \frac{\partial }{\partial x_i}
\left(\frac{x_k}{\rho^\alpha}\right)\, dx
+\int_\OO \frac{\partial^2 u}{\partial x_i\partial x_j}\cdot
\frac{\partial u}{\partial x_k}\cdot
\frac{\partial^2}{\partial x_i\partial x_j}
\left(\frac{x_k}{\rho}\right)\, dx.
\endaligned
\tag 3.12
$$
Substituting (3.10) and (3.12) into (3.9) and using
$$
\aligned
&\frac{\partial }{\partial x_i}
\left(\frac{x_j}{\rho^\alpha}\right)
=\frac{\delta_{ij}}{\rho^\alpha}
-\alpha \, \frac{ x_i x_j}{\rho^{\alpha +2}},\\
&\frac{\partial^2}{\partial x_i\partial x_j}
\left(\frac{x_k}{\rho^\alpha}\right)
=-\alpha (\delta_{ik} x_j +\delta_{ij} x_k +\delta_{jk} x_i)\rho^{-\alpha -2}
+\alpha (\alpha +2)\, \frac{x_i x_j x_k}{\rho^{\alpha +4}},
\endaligned
\tag 3.13
$$
we obtain
$$
\aligned
&\int_\OO \Delta^2 u\cdot \frac{\partial u}{\partial\rho}
\, \frac{dx}{\rho^{\alpha -1}}\\
&\ \ \ \ \ \ 
=-\frac12 \int_{\partial\OO} |\nabla^2 u|^2\, <\frac{x}{\rho}, N>
\, \frac{d\sigma}{\rho^{\alpha -1}}
+\frac12 (\alpha +4 -n)\int_\OO |\nabla^2 u|^2\, \frac{dx}{\rho^\alpha}\\
&\ \ \ \ \ \ \ \ \ \ \   
-2\alpha \int_\OO \big|\frac{\partial }{\partial\rho}
\nabla u\big|^2\, \frac{dx}{\rho^\alpha}
-2\alpha \int_\OO \frac{\partial^2 u}{\partial x_i\partial x_j}
\cdot \frac{\partial u}{\partial x_j}\cdot\frac{x_i}{\rho^{\alpha +2}}\,
dx\\
&\ \ \ \ \ \ \ \ \ \ \ \ 
-\alpha \int_\OO \Delta u\cdot \frac{\partial u}{\partial \rho}\,
\frac{dx}{\rho^{\alpha +1}}
+\alpha (\alpha +2)\int_\OO
\frac{\partial^2 u}{\partial x_i\partial x_j}\cdot
\frac{\partial u}{\partial x_k}\cdot
\frac{x_i x_j x_k}{\rho^{\alpha +4}}\, dx.
\endaligned
\tag 3.14
$$

Finally we note that
$$
\aligned
 2\int_\OO \frac{\partial^2 u}{\partial x_i\partial x_j}
\cdot \frac{\partial u}{\partial x_j}\cdot
\frac{x_i}{\rho^{\alpha +2}}\, dx
&=-\int_\OO |\nabla u|^2\, \frac{\partial}{\partial x_i}
\left(\frac{x_i}{\rho^{\alpha +2}}\right)\, dx\\
&
= (\alpha +2-n)\int_\OO |\nabla u|^2\, \frac{dx}{\rho^{\alpha +2}},
\endaligned
\tag 3.15
$$
and
$$
\aligned
\int_\OO \frac{\partial^2 u}{\partial x_i
\partial x_j}\cdot \frac{\partial u}{\partial x_k}
\cdot \frac{x_i x_j x_k}{\rho^{\alpha +4}}\, dx
&=-\frac12 
\int_\OO \frac{\partial u}{\partial x_j}
\cdot\frac{\partial u}{\partial x_k}
\cdot \frac{\partial }{\partial x_i}
\left(\frac{x_i x_j x_k}{\rho^{\alpha +4}}\right)\, dx\\
&=\frac12
(\alpha +2 -n)\int_\OO \big|\frac{\partial u}
{\partial\rho}\big|^2\, \frac{dx}{\rho^{\alpha +2}}.
\endaligned
\tag 3.16
$$
The desired formula (3.8) follows by substituting (3.15), (3.16) as well as
(2.19) into (3.14).
The proof is complete.
\enddemo

\proclaim{\bf Lemma 3.17}
Suppose that $u\in C^1(\overline{\OO})$ and $u=0$ on $\partial\OO$.
Then, for any $\alpha\in \br$,
$$
\int_\OO \big|\frac{\partial}{\partial \rho}
\left(u\rho^{\frac{n-\alpha}{2}}\right)\big|^2\, \frac{dx}{\rho^{n-2}}
=\int_\OO \big|\frac{\partial u}{\partial\rho}\big|^2\, 
\frac{dx}{\rho^{\alpha -2}}
-\frac14 (n-\alpha)^2 \int_\OO |u|^2\, \frac{dx}{\rho^\alpha},
\tag 3.18
$$
where $\rho=|x-y|$ with $y\in \overline{\OO}^c$ fixed.
\endproclaim

\demo{Proof} To see (3.18), we note that
$$
\aligned
&\big|\frac{\partial}{\partial \rho}
\left(u\rho^{\frac{n-\alpha}{2}}\right)\big|^2\\
&=\big|\frac{\partial u}{\partial\rho}\big|^2
\,\rho^{n-\alpha}
+(n-\alpha)\, u\, \frac{\partial u}{\partial\rho}\, \rho^{n-\alpha-1}
+\frac14 (n-\alpha)^2\, |u|^2\, \rho^{n-\alpha -2}.
\endaligned
\tag 3.19
$$
Also, using integration by parts and $u=0$ on $\partial\OO$, we have
$$
\aligned
(n-\alpha)\int_\OO u\, \frac{\partial u}{\partial\rho}\,
\rho^{n-\alpha -1}\, \frac{dx}{\rho^{n-2}}
&=-\frac12 (n-\alpha)\int_\OO |u|^2\, \frac{\partial}
{\partial x_i}\left(\frac{x_i}{\rho^\alpha}\right)\, dx\\
&=-\frac12 (n-\alpha)^2
\int_\OO |u|^2\, \frac{dx}{\rho^\alpha}.
\endaligned
\tag 3.20
$$
In view of (3.19), this gives (3.18).
\enddemo

We are now ready to prove the integral identity (3.1).

\proclaim{\bf Lemma 3.21}
Let $\OO$ be a bounded Lipschitz domain in $\br^n$, $n\ge 2$.
Suppose that $u\in C^4(\overline{\OO})$ and
$u=0$, $\nabla u=0$ on $\partial\OO$. Then (3.1) holds for
any $\alpha<n$ and any $y\in \partial\OO$.
\endproclaim

\demo{Proof} By the Lebesgue Dominated Convergence Theorem, it suffices
to establish (3.1) for $y\in \overline{\OO}^c$.
To this end, we note that
$$
\int_\OO \Delta u \cdot \Delta \left(\frac{u}{\rho^\alpha}\right)\, dx
=\int_\OO \frac{\partial^2 u}{\partial x_i \partial x_j}
\cdot
\frac{\partial^2 u}
{\partial x_i\partial x_j} \left(\frac{u}{\rho^\alpha}\right)\, dx,
\tag 3.22
$$
from integration by parts. Thus, by (3.3) and (3.8), we have
$$
\aligned
(\alpha +4-n)\int_\OO \Delta u\cdot &\Delta \left(\frac{u}{\rho^\alpha}
\right)\, dx
-2\int_\OO \Delta^2 u\cdot \frac{\partial u}{\partial \rho}\, 
\frac{dx}{\rho^{\alpha-1}}\\
&=
\int_{\partial\OO} |\nabla^2 u|^2 <\frac{x-y}{\rho},N>\,
\frac{d\sigma}{\rho^{\alpha -1}}
\\
& \ \ \  
+4\alpha \int_\OO \big|\frac{\partial }{\partial\rho}
\nabla u\big|^2\, \frac{dx}{\rho^\alpha}
-\alpha (n-\alpha -2)^2\int_\OO |\nabla u|^2\, \frac{dx}{
\rho^{\alpha +2}}\\
&\ \ \ \ 
+2\alpha (\alpha+2)(n-\alpha-2)\int_\OO \big|
\frac{\partial u}{\partial\rho}\big|^2\, \frac{dx}{\rho^{\alpha +2}}\\
&\ \ \ \ \ -\frac12
\alpha (\alpha+2)(n-\alpha-2)(n-\alpha -4)^2
\int_\OO |u|^2\, \frac{dx}{\rho^{\alpha +4}}.
\endaligned
$$
In view of (3.18), this gives the integral identity (3.1).
\enddemo

Next we will use (3.1) to derive estimate (1.4) on convex domains
with smooth boundaries for any $\lambda<n$.

\proclaim{\bf Lemma 3.23}
Let $\OO$ be a convex domain in $\br^n$, $n\ge 4$ with smooth boundary.
Let $0<\lambda<n$. Then there exist constants $C_0>0$ and $R_0>0$
depending only on $n$, $\lambda$ and the Lipschitz character of
$\OO$ such that
estimate (1.4) holds for any $v$ satisfying (1.5).
\endproclaim

\demo{Proof} Let $R_0>0$ be a constant so that for any $Q\in \partial
\OO$, $T(Q,4R_0)$ is given by  
the intersection of $B(Q,4R_0)$ and the region above a 
Lipschitz graph, after a possible rotation.
Fix $Q\in \partial\OO$ and $0<R<R_0$. Let $v$ be a biharmonic function in
$\OO$ such that $v=\frac{\partial v}{\partial N}=0$ on $I(Q,R)$ and
$(\nabla v)^*\in L^2(\partial\OO)$. Since $\OO$ has smooth boundary,
by the classical regularity theory for elliptic equations,
$v\in C^4(\overline{T(Q,R/2)})$.

Let $\eta$ be a smooth function on $\br^n$ such that
$\eta=1$ on $B(Q,R/8)$, supp$\, \eta\subset B(Q,R/4)$ and
$|\nabla^k \eta|\le C/R^k$ for $0\le k\le 4$.
 Note that $u=v\eta\in C^4(\overline{\OO})$
and $u=0$, $\nabla u=0$ on $\partial\OO$. Thus we may apply
integral identity (3.1) to $u$ with $\alpha =n-2$
and $y=Q$. This gives
$$
\int_\OO \big|\frac{\partial}{\partial\rho}
(\nabla u)\big|^2\, \frac{dx}{\rho^{n-2}}
\le C\, \left\{
\left|\int_\OO \Delta^4 u\cdot \frac{\partial u}{\partial \rho}\,
\frac{dx}{\rho^{n-3}}\right|
+\left|\int_\OO \Delta^4 u\cdot u \, \frac{dx}{\rho^{n-2}}\right|\right\}.
\tag 3.24
$$
Since $\Delta^2 v=0$ in $\OO$, we have
$$
\Delta^2 u=2<\nabla(\Delta v),\nabla\eta>
+\Delta v\cdot \Delta \eta
+\Delta\left\{ 2<\nabla v,\nabla\eta> +v\Delta\eta\right\}.
\tag 3.25
$$
Substituting (3.25) into (3.24) and using integration by parts
as well as Cauchy inequality,
we obtain
$$
\aligned
\int_\OO \big|\frac{\partial}{\partial\rho}
(\nabla u)\big|^2\, \frac{dx}{\rho^{n-2}}
&\le \frac{C}{R^{n-2}}
\int_{T(Q,R/4)}\left\{ |\nabla^2 v|^2 +\frac{|\nabla v|^2}{R^2}
+\frac{|v|^2}{R^4}\right\}\, dx\\
&\le \frac{C}{R^{n}}
\int_{T(Q,R/2)}|\nabla v|^2\, dx,
\endaligned
\tag 3.26
$$
where we also used the Cacciopoli's inequality (2.2) and
Poincar\'e inequality in the second inequality.

Since supp$(u)\subset B(Q,R)$, for any $\delta\in (0,n-2)$, we have
$$
\int_\OO \big|\frac{\partial}{\partial\rho}
(\nabla u)\big|^2\, \frac{dx}{\rho^{n-2}}
\ge \frac{1}{R^\delta}
\int_\OO \big|\frac{\partial}{\partial\rho}
(\nabla u)\big|^2\, \frac{dx}{\rho^{n-2-\delta}}
\ge \frac{\delta^2}{4R^\delta}
\int_\OO |\nabla u|^2\, \frac{dx}{\rho^{n-\delta}},
$$
where the second inequality follows from (3.18) with $\alpha
=n+2-\delta$, which also holds
for $y\in \partial\OO$ if $\alpha<n+2$.
In view of (3.26), this gives
$$
\int_{T(Q,r)} |\nabla v|^2\, dx
\le r^{n-\delta} \int_{T(Q,r)}|\nabla u|^2\, \frac{dx}{\rho^{n-\delta}}
\le C_\delta \, \left(\frac{r}{R}\right)^{n-\delta}
\int_{T(Q,R)} |\nabla v|^2\, dx,
$$
for any $0<r<R/8$. Estimate (1.4) is thus proved for
$\lambda=n-\delta$. 
\enddemo

Lemma 3.23, together with a well known approximation argument, gives
part (b) of the Main Theorem.

\proclaim{\bf Theorem 3.27}
Let $\OO$ be a bounded convex domain in $\br^n$, $n\ge 4$.
Then the $L^p$ Dirichlet problem (1.1) is uniquely solvable
for any $2-\e<p<\infty$.
\endproclaim

\demo{Proof}
Let $p>2$ and $f\in L^p_1(\partial\OO)$, $g\in L^p(\partial\OO)$. We need
to show that the unique solution $u$ to the $L^2$ Dirichlet problem
(1.1) satisfies 
estimate (1.7). To this end, we first note that
by an approximation argument(e.g. see \cite{JK} for
Laplace's equation), we may assume that
$f,g\in C_0^\infty(\br^n)$.

Next we approximate $\OO$ from outside by a sequence of convex domains
$\{ \OO_j\}$ with smooth boundaries,
$\OO_1\supset \OO_2\supset \dots \supset\OO$.
Let $u_j$ be the solution to the $L^2$ Dirichlet problem
(1.1) on $\OO_j$ with boundary data $(u_j,\frac{\partial u_j}{
\partial N})=(f|_{\partial\OO_j},
g|_{\partial\OO_j})$ on $\partial\OO_j$.
 By Lemma 3.23 and Theorem 1.3, we have
$$
\| (\nabla u_j)^*\|_{L^p(\partial\OO)}
\le C\, \| (\nabla u_j)^*_j\|_{L^p(\partial\OO_j)}
\le C\, \left\{ \|\nabla_t f\|_{L^2(\partial\OO_j)}
+\| g\|_{L^p(\partial\OO_j)}\right\},
\tag 3.28
$$
where $(\nabla u_j)^*_j$ denotes the non-tangential maximal function
of $\nabla u_j$ with respect to $\OO_j$, and $C$ is a constant independent
of $j$.
Estimate (3.28) implies
 that the sequence $\{ \nabla u_j\}$ is uniformly bounded on any compact
subset of $\OO$. It follows that there exist
 a subsequence, which we still
denoted by $\{ \nabla u_j\}$, and a function $u$ on $\OO$ such that
$u_j$ converges to $u$ uniformly on any compact subset of
$\OO$. 
It is easy to show that $u$ is biharmonic in $\OO$. 
Also by (3.28) and Fatou's Lemma,
$$
\| (\nabla u)^*_K\|_{L^p(\partial\OO)}
\le C\, \left\{
\|\nabla_t f\|_{L^p(\partial\OO)}
+\| g\|_{L^p(\partial\OO)}\right\},
\tag 3.29
$$
where $K$ is a compact subset of $\OO$, and 
$(\nabla u)^*_K(Q)=\sup\{ |\nabla u(x)|:\
x\in K \text{ and } |x-Q|< 2\,
\text{dist}(x,\partial\OO)\}$.
By the monotone convergence theorem, this gives
 the estimate (1.7) on $\OO$.

Finally one may use $L^2$ estimates on $\| (\nabla u_i -\nabla u_j)^*
\|_{L^2(\partial\OO_j)}$ for $i\ge j$ as well as
$L^2$
regularity estimate, $\| (\nabla^2 u_j)_j^*\|_{L^2(\partial\OO_j)}
\le C\, \{ \|\nabla^2 f\|_{L^2(\partial\OO_j)}
+\|\nabla g\|_{L^2(\partial\OO_i)}\}$
(see \cite{V2})
to show that $u=f$ and $\frac{\partial u}{\partial N}=g$
on $\partial \OO$ in the sense of non-tangential convergence.
We leave the details to the reader.
\enddemo
\enddemo

\enddocument

\end